\documentclass[a4paper,10pt,draft]{amsart}

\usepackage{verbatim, amssymb, enumerate}

\renewcommand \a{\alpha}
\renewcommand \b{\beta}
\newcommand \K{\delta}
\newcommand \n{\nabla}
\newcommand \la{\lambda}
\newcommand \ve{\varepsilon}
\newcommand \id{\mathrm{id}}
\newcommand \br{\mathbb{R}}
\newcommand \bc{\mathbb{C}}

\newcommand \rk{\operatorname{rk}}
\newcommand \Ker{\operatorname{Ker}}

\newcommand \End{\operatorname{End}}

\newcommand \Span{\operatorname{Span}}
\newcommand \Tr{\operatorname{Tr}}
\newcommand \db{\partial}

\newcommand \cp{\mathcal{C}}
\newcommand \og{\mathfrak{o}}
\newcommand\ag{\mathfrak a}
\newcommand\g{\mathfrak g}
\newcommand\h{\mathfrak h}

\renewcommand\t{\mathfrak t}
\newcommand \Sk{\mathrm{Skew}}

\newcommand \Hom{\operatorname{Hom}}

\newcommand \ad{\operatorname{ad}}

\newcommand \diag{\operatorname{diag}}

\renewcommand \O{\operatorname{O}}

\newcommand \<{\langle}
\renewcommand \>{\rangle}

\newcommand \mU{\mathcal{U}}

\theoremstyle{plane}

\newtheorem*{theorem*}{Theorem}
\newtheorem*{corollary*}{Corollary}
\newtheorem*{conj*}{Conjecture}
\newtheorem{lemma}{Lemma}

\newtheorem*{prop*}{Proposition}

\theoremstyle{definition}

\newtheorem*{definition*}{Definition}

\theoremstyle{remark}

\begin{document}

\title{Weyl homogeneous manifolds modeled on compact Lie groups}

\author{Y.Nikolayevsky}
\address{Department of Mathematics and Statistics, La Trobe University, Victoria, 3086, Australia}
\email{y.nikolayevsky@latrobe.edu.au}

\date{\today}


\subjclass[2000]{Primary: 53C35,53A30, secondary: 53B20}
\keywords{Weyl tensor, symmetric space}

\begin{abstract}
A Riemannian manifold is called Weyl homogeneous, if its Weyl tensors at any two points are ``the same", up to
a positive multiple. A Weyl homogeneous manifold is modeled on a homogeneous space $M_0$, if its Weyl tensor at
every point is ``the same" as the Weyl tensor of $M_0$, up to a positive multiple.
We prove that a Weyl homogeneous manifold $M^n, \; n \ge 4$, modeled on an irreducible symmetric space
$M_0$ of types II or IV (compact simple Lie group with a bi-invariant metric or its noncompact dual) is conformally
equivalent to $M_0$.
\end{abstract}

\maketitle

\section{Introduction}
\label{s:intro}

The question which we discuss in this paper is inspired firstly, by an attempt to generalize the classical Weyl-Schouten Theorem to a wider class of model
spaces, and secondly, by translating the notion of curvature homogeneity to the Weyl tensor. Recall that a smooth Riemannian manifold $M^n$ is called \emph{curvature homogeneous}, if for
any $x, y \in M^n$, there exists a linear isometry $\iota:T_xM^n \to T_yM^n$ which maps the curvature tensor of $M^n$ at $x$ to the curvature tensor of
$M^n$ at $y$. A smooth Riemannian manifold $M^n$ is \emph{modeled on a homogeneous space $M_0$}, if for every point $x \in M^n$, there exists a linear
isometry $\iota:T_xM^n \to T_oM_0$ which maps the curvature tensor of $M^n$ at $x$ to the curvature tensor of $M_0$ at $o \in M_0$ (such a manifold is
then automatically curvature homogeneous). The term ``curvature homogeneous" was coined by F.Tricerri and L.Vanhecke in 1986 \cite{TV}.
Since then, the theory of curvature homogeneous manifolds received a great deal of attention of geometers; for the current state of knowledge in the area the reader is
referred to the excellent book \cite{Gil}. We adopt the following definition.

\begin{definition*} \label{d:ws}
A smooth Riemannian manifold $M^n$ is called \emph{Weyl homogeneous}, if for any $x, y \in M^n$, there exists a linear isometry
$\iota:T_xM^n \to T_yM^n$ which maps the Weyl tensor of $M^n$ at $x$ to a positive multiple of the Weyl tensor of $M^n$ at $y$.
A smooth Weyl homogeneous Riemannian manifold $M^n$ is \emph{modeled on a homogeneous space $M_0$}, if for every point
$x \in M^n$, there exists a linear isometry $\iota:T_xM^n \to T_oM_0$ which maps the Weyl tensor of $M^n$ at $x$ to a
positive multiple of the Weyl tensor of $M_0$ at $o \in M_0$.
\end{definition*}

In the latter case, we will also say that $M^n$ has \emph{the same} Weyl tensor as $M_0$.

Clearly, a Riemannian manifold which is conformally equivalent to a homogeneous space is Weyl homogeneous. The converse is
not true even for Weyl homogeneous manifolds modeled on symmetric spaces (see Section~\ref{s:red}, where we discuss an example
from \cite[Theorem~4.2]{BKV} from the conformal point of view). Moreover, based on the existence of many examples of curvature
homogeneous manifolds, which are not locally homogeneous \cite{Gil}, one should most probably expect at least as many examples
in the conformal settings.

We will therefore restrict ourselves to the Weyl homogeneous manifolds modeled on symmetric spaces and will consider the following question:
\begin{quote}
\emph{For which symmetric spaces $M_0$, a Riemannian manifold having the same Weyl conformal curvature tensor as $M_0$,
is conformally equivalent to $M_0$?}
\end{quote}
The answer to this question is in positive, when $M_0=\br^n, \, n \ge 4$ (or more generally, when $M_0$ is a conformally flat symmetric space
of dimension $n \ge 4$) and when $M_0$ is a rank-one symmetric space of dimension $n > 4$ \cite[Theorem~2]{N}.
Very likely candidates for a positive answer would also be irreducible symmetric spaces, as by \cite{TV}, a curvature homogeneous Riemannian
manifold modeled on an irreducible symmetric space, is locally isometric to it. In this paper, we prove that this is indeed the case when the model space
$M_0$ is an irreducible symmetric space of types II or IV, that is, a compact simple Lie group with a bi-invariant metric or its noncompact dual:

\begin{theorem*}
Let $M^n, \; n \ge 4$, be a smooth Weyl homogeneous Riemannian manifold modeled on an irreducible symmetric space $M_0$
of types II or IV. Then $M^n$ is locally conformally equivalent to $M_0$.
\end{theorem*}
The requirement $n \ge 4$ is, of course, necessary. Note that for arbitrary irreducible symmetric spaces, one should require at least $n >4$ for the
claim of the Theorem to remain true. This is because a four-dimensional Riemannian manifold having the same Weyl tensor as $\bc P^2$ is either self-dual
or anti-self-dual by \cite{BG2} and because there exist self-dual K\"{a}hler metrics on $\bc^2$ which are not locally conformally equivalent to
locally symmetric ones \cite{Der}.

\smallskip

The paper is organized as follows.
In Section~\ref{s:geo}, after a brief introduction, we prove the Theorem with the help of a purely Lie-algebraic Proposition whose
proof, in turn, is given in Section~\ref{s:alg}. In Section~\ref{s:red}, we verify that the claim of the Theorem is false, if $M_0=M^2(\kappa) \times \br^{n-2}$,
where $M^2(\kappa)$ is the two-dimensional space of the constant curvature $\kappa \ne 0$ (this should not be surprising, as there exist curvature
homogeneous Riemannian spaces modeled on such an $M_0$ which are not locally symmetric \cite[Theorem~4.2]{BKV}).

The Riemannian manifold $M^n$ is assumed to be smooth (of class $C^\infty$), although the Theorem remains valid for
manifolds of class $C^k$, with sufficiently large $k$. 

\section{Proof of the Theorem}
\label{s:geo}

Let $M^n$ be a Riemannian manifold with the metric $\<\cdot, \cdot\>$ and the Levi-Civita connection $\n$. For
vector fields $X,Y$, define the field of linear operators $X \wedge Y$ by lowering the index of the corresponding
bivector: $(X \wedge Y)Z=\<X,Z\>Y-\<Y,Z\>X$. The curvature tensor is defined by $R(X,Y)=\n_X\n_Y-\n_Y\n_X-\n_{[X,Y]}$,
where $[X,Y]=\n_XY-\n_YX$, and the Weyl conformal curvature tensor $W$, by
\begin{equation}\label{eq:weyldef}
R(X,Y)=(\rho X) \wedge Y + X \wedge (\rho Y) + W(X, Y), 
\end{equation}
where $\rho =\frac{1}{n-2}\operatorname{Ric}-\frac{\operatorname{scal}}{2(n-1)(n-2)}\id$,
$\operatorname{Ric}$ is the Ricci operator and $\operatorname{scal}$ is the scalar curvature. We denote
$W(X,Y,Z,V)=\<W(X,Y)Z,V\>$.

We start with the following general fact.
\begin{lemma}\label{l:eismooth}
Suppose that $M^n$ is a Weyl homogeneous manifold with the metric $\<\cdot, \cdot\>'$ modeled on a homogeneous space $M_0$
with the Weyl tensor $W_0 \ne 0$. Choose a point $o \in M_0$ and an orthonormal basis $E_i$ for $T_oM_0$. 
Then there exists a smooth metric $\<\cdot, \cdot\>$ on $M^n$ conformally equivalent to $\<\cdot, \cdot\>'$ such that for every $x \in M^n$,
there exists a smooth orthonormal frame $e_i$ (relative to $\<\cdot, \cdot\>$) on a neighborhood $\mU=\mU(x) \subset M^n$ satisfying
$W(e_i,e_j,e_k,e_l)(y)=W_0(E_i,E_j,E_k,E_l)$, for all $y \in \mU$.
\end{lemma}
\begin{proof}
By definition, there exists a function $f: M^n \to \br^+$ and, at every point $x \in M^n$, there exist an orthonormal
basis $e'_i$ for $T_xM^n$ (relative to $\<\cdot, \cdot\>'$) such that
$W'(e'_i,e'_j,e'_k,e'_l)=f(x)W_0(E_i,E_j,E_k,E_l)$, for all $i,j,k,l=1, \dots, n$, where $W'$ is the Weyl tensor for $\<\cdot, \cdot\>'$.
It follows that the function $f= \|W'\|'\|W_0\|^{-1}$ is smooth, so the metric $\<\cdot, \cdot\>= f^{-1/2}\<\cdot, \cdot\>'$ is smooth, and
moreover, for the vectors $e_i= f^{1/4}e'_i$ (which are orthonormal relative to $\<\cdot, \cdot\>$), we have
$W(e_i,e_j,e_k,e_l)=W_0(E_i,E_j,E_k,E_l)$, where $W=W'$ is the Weyl tensor for $\<\cdot, \cdot\>$. In the remaining part of the proof
we assume that the metric on $M^n$ is $\<\cdot, \cdot\>$.

Define the map $\tau: \O(n) \to \br^{n^4}$ by $\tau(g)_{ijkl}=W_0(gE_i,gE_j,gE_k,gE_l)$, where $1 \le i,j,k,l \le n$ and $g \in \O(n)$.
The group $H=\{g \in \O(n) \, : \, \tau(g)=\tau(\id)\}$ is a compact subgroup of $\O(n)$, and $\tau(\O(n))$, the orbit 
of $W_0$, is diffeomorphic to the homogeneous space $\O(n)/H$. 
Choosing a smooth orthonormal frame $e_i$ on a small neighborhood $\mU=\mU(x)$, we obtain a smooth map $F: \mU \to \br^{n^4}$
defined by $F(y)_{ijkl}=W(e_i,e_j,e_k,e_l)$, for $y \in \mU$. The image of $F$ lies in $\tau(\O(n))=\O(n)/H$, which gives a smooth map from $\mU$ to $\O(n)/H$ that
can be locally lifted to a smooth map $\phi: \mU \to \O(n)$,
so that $W(e_i,e_j,e_k,e_l)(y)=W_0(\phi(y)E_i,\phi(y)E_j,\phi(y)E_k,\phi(y)E_l)$, for $y \in \mU$. Then $\phi(y)^{-1}e_i$ is
the required smooth orthonormal frame on $\mU$.
\end{proof}

For the remainder of the proof we will assume that the metric on $M^n$ is chosen as in Lemma~\ref{l:eismooth} and we will be proving that
$M^n$ (with that metric) is locally isometric to the model space $M_0$.

Let $x \in M^n$ and let $e_i$ be the orthonormal frame on the neighborhood $\mU$ of $x$ introduced in Lemma~\ref{l:eismooth}.
For every $Z \in T_xM^n$, define the linear operator $K_Z$ (the connection operator) by
\begin{equation}\label{eq:gdefK}
K_Ze_i=\n_Ze_i
\end{equation}
(and extended to $T_xM^n$ by linearity). As the basis $e_i$ is orthonormal, $K_Z$ is skew-symmetric, for all $Z \in T_xM^n$.
For smooth vector fields $X, Y$ on $\mU$ define
\begin{equation}\label{eq:gdefPhi}
\Phi(X,Y)=(\n_X \rho)Y-(\n_Y \rho)X,
\end{equation}
where $\rho =\frac{1}{n-2}\operatorname{Ric}-\frac{\operatorname{scal}}{2(n-1)(n-2)}\id$ (see \eqref{eq:weyldef}).
As $\rho$ is symmetric, we also have
\begin{equation}\label{eq:cyclePhi}
\sigma_{XYZ}\<\Phi(X,Y),Z\>=0,
\end{equation}
where $\sigma_{XYZ}$ is the sum over the cyclic permutations of $X,Y,Z$.

Let $M_0=\mathbf{G}/\mathbf{H}$ be the model (simply connected) symmetric space for $M^n$, where $\mathbf{G}$ is the full isometry group of $M_0$ and $\mathbf{H}$ is the
isotropy subgroup of $o \in M_0$, and let $\mathfrak{G}=\h+\mathfrak{t}$ be the corresponding Cartan
decomposition, where $\mathfrak{G}$ and $\h$ are the Lie algebras of $\mathbf{G}$ and $\mathbf{H}$
respectively, and $\mathfrak{t}=T_oM_0$. Denote $R_0$ the curvature tensor of $M_0$ at $o$, so that for
$X, Y, Z \in T_oM_0, \; R_0(X,Y)Z=-[[X,Y],Z]=-\ad^{\h}_{[X,Y]}Z$.

In the assumptions of Lemma~\ref{l:eismooth}, identify $T_xM^n$ with $T_oM_0$ via the linear isometry
$\iota$ mapping $e_i$ to $E_i$. Define $K$ and $\Phi$ on $\t=T_oM_0$, by the pull-back by $\iota$.

Let $\operatorname{Ric_0}$ and $\operatorname{scal}_0$ be the Ricci tensor and the scalar curvature
of $M_0$ (at $o \in M_0$), and let $\rho_0=\frac{1}{n-2}\operatorname{Ric_0}-\frac{\operatorname{scal}_0}{2(n-1)(n-2)}\id$ (see \eqref{eq:weyldef}).
Define the operator $\Psi: \Lambda^2\t \to \t$ by
\begin{equation}\label{eq:defPsi}
\Psi(X,Y)=\Phi(X,Y) + [\rho_0,K_X]Y - [\rho_0,K_Y]X,
\end{equation}
where (here and below) the bracket of linear operators is the usual commutator. From \eqref{eq:cyclePhi} and the fact that $[\rho_0,K_X]$ is symmetric it follows that
\begin{equation}\label{eq:cyclePsi}
\sigma_{XYZ}\<\Psi(X,Y),Z\>=0, \quad \text{for all} \quad  X,Y,Z \in \t.
\end{equation}

\begin{lemma}\label{l:gbianchi}
In the assumptions of Lemma~\ref{l:eismooth}, let $M_0$ be a symmetric space. For $x \in M^n$, identify $T_xM^n$ with $\t=T_oM_0$ via the linear isometry
$\iota$ mapping $e_i$ to $E_i$. Define $K$ and $\Phi$ on $T_xM^n$ by (\ref{eq:gdefK}, \ref{eq:gdefPhi}) and on $T_oM_0$, by
the pull-back by $\iota$, and define $\Psi$ by \eqref{eq:defPsi}. Then
\begin{gather}\label{eq:gbiad}
\sigma_{XYZ}([\ad^{\h}_{[X, Y]},K_Z]+\ad^{\h}_{[K_XY-K_YX, Z]}+ \Psi(X,Y)\wedge Z)=0, \\
(\n_Z W)(X,Y)=[\ad^{\h}_{[X, Y]},K_Z]+\ad^{\h}_{[K_XY-K_YX, Z]}+([\rho_0,K_Z]X) \wedge Y - X \wedge ([\rho_0,K_Z]Y). \label{eq:gnablaW}
\end{gather}
\end{lemma}

\begin{proof}
For the orthonormal frame $e_i$ on $\mU$ introduced in Lemma~\ref{l:eismooth}, $W(e_i,e_j,e_k,e_l)=W_0(E_i,E_j,E_k,E_l)$,
which implies $(\n_Z W)(e_i,e_j,e_k,e_l)=-W(K_Ze_i,e_j,e_k,e_l)- W(e_i,K_Ze_j,e_k,e_l)-W(e_i,e_j,K_Ze_k,e_l)$ $-W(e_i,e_j,e_k,K_Ze_l)$.
Then for $X,Y,Z \in T_xM^n$,
\begin{equation}\label{eq:nWonM}
(\n_Z W)(X,Y)=-W(K_ZX,Y)- W(X,K_ZY)-[W(X,Y),K_Z],
\end{equation}
so by \eqref{eq:weyldef}, $(\n_ZR)(X,Y)=(\n_Z\rho X)\wedge Y + X \wedge (\n_Z\rho Y) - W(K_ZX, Y)-W(X, K_ZY)-[W(X, Y),K_Z]$. Then by
\eqref{eq:gdefPhi}, the second Bianchi identity gives
\begin{equation}\label{eq:gbionM}
\sigma_{XYZ}(\Phi(X, Y)\wedge Z - W(K_ZX, Y)-W(X, K_ZY)-[W(X, Y),K_Z])=0.
\end{equation}
Identifying $T_xM^n$ with $T_oM_0$ via $\iota$ we have $W(X, Y)Z=W_0(X, Y)Z$, for $X,Y,Z \in T_xM^n$,
where $W_0$ is the Weyl tensor of $M_0$ at $o$. Then by \eqref{eq:weyldef},
$W(X, Y)=R_0(X, Y)-(\rho_0 X) \wedge Y - X \wedge (\rho_0 Y)=-\ad^{\h}_{[X, Y]}-(\rho_0 X) \wedge Y - X \wedge (\rho_0 Y)$.
Substituting this to \eqref{eq:gbionM}
(and to \eqref{eq:nWonM}) and using \eqref{eq:defPsi} and the fact that $[K_Z,X \wedge Y]=(K_ZX) \wedge Y+X \wedge (K_ZY)$
we obtain \eqref{eq:gbiad} (and \eqref{eq:gnablaW}).
\end{proof}

For a Euclidean space $V$, denote $\Sk(V)=\og(V)$ the space of the skew-symmetric operators on $V$, which we will identify with the
space $\Lambda^2 V$ of bivectors using the inner product. For a Lie algebra $\g$, $\ad(\g) \subset \End(\g)$ is the space of derivations of $\g$
(if the inner product is bi-invariant,  $\ad(\g) \subset \Sk(\g)$).

The next step in the proof requires the following algebraic proposition whose prove will be given in Section~\ref{s:alg}.

\begin{prop*}
Let $\g$ be a simple compact Lie algebra of dimension $n > 4$ with the inner product being the negative of the Killing
form. Suppose that the maps
$K \in \Hom(\g, \Sk(\g)), \; K: Z \to K_Z$ and $\Phi \in \Hom(\Lambda^2\g, \g), \; \Phi: X \wedge Y \to \Phi(X,Y)$
satisfy \eqref{eq:cyclePhi} and
\begin{equation}\label{eq:biad}
\sigma_{XYZ}([\ad_{[X, Y]},K_Z]+\ad_{[K_XY-K_YX, Z]}+ \Phi(X,Y)\wedge Z)=0,
\end{equation}
for all $X,Y,Z \in \g$. Then $\Phi=0$ and $K \in \Hom(\g, \ad(\g))$.
\end{prop*}

With the Proposition, we can finish the proof of the Theorem as follows. Suppose that $M_0$ is an irreducible symmetric space of types II or IV.
In the both cases, as $M_0$ is Einstein, we can take $\Psi=\Phi$ in \eqref{eq:gbiad} and drop the last two terms in \eqref{eq:gnablaW}.

If $M_0$ is an irreducible symmetric space of type II (a compact simple Lie group $G$ with a biinvariant metric), then
$\mathfrak{G}=\g \oplus \g$, where $\g$ is the Lie algebra of $G$ and $\h = \g$ is the diagonal in $\mathfrak{G}$
and $T_oM_0=\mathfrak{t}=\{(T,-T)\, : \, T \in \g\}=\g$ (as a linear space). Identifying $T_oM_0$ with $\g$ by projecting to
the  first component we get $\ad^{\h}_{[X,Y]}Z=\ad_{[X,Y]}Z$, for $X,Y,Z \in T_oM_0$. Then equation \eqref{eq:gbiad} becomes
equation \eqref{eq:biad}, so by the Proposition, $\Phi=0$ and $K_Z=\ad^{\h}_{AZ}$, for some linear map $A: \mathfrak{t} \to \h$, so $\n W=0$ by \eqref{eq:gnablaW}.
If $M_0$ is an irreducible symmetric space of type IV, then $\mathfrak{G}=\g \oplus \mathrm{i}\g$, where $\g$ is the Lie algebra of a compact simple Lie group $G$,
$\h = \g$ and $T_oM_0=\mathfrak{t}=\mathrm{i}\g$. Identifying $T_oM_0$ with $\g$ we get $\ad^{\h}_{[X,Y]}Z=-\ad_{[X,Y]}Z$, for $X,Y,Z \in T_oM_0$. Then
\eqref{eq:gbiad} becomes \eqref{eq:biad}, if we replace $\Phi$ by $-\Phi$. So by the Proposition, $\Phi=0$ and $K_Z=\ad^{\h}_{AZ}$,
for some $A: \mathfrak{t} \to \h$, which again implies $\n W=0$ by \eqref{eq:gnablaW}.

By \cite{Rot}, as $\n W=0$, but $W \ne 0$, the manifold $M^n$ is locally symmetric. To prove that $M^n$ is locally isometric to
$M_0$, note that the equation $\Phi=0$ and \eqref{eq:gdefPhi} imply that $\rho$ is a (symmetric) Codazzi tensor. Then by
\cite[Theorem~1]{DS}, the exterior products of the eigenspaces of $\rho$ are invariant subspaces of the curvature operator on
the bivectors. Then by \eqref{eq:weyldef}, they are also invariant subspaces of the Weyl operator, hence (by \eqref{eq:weyldef}
again and by the fact that $W=W_0$ and $M_0$ is Einstein) they are invariant subspaces of the (pull-back of the) curvature operator
$R_0$ on the bivectors. Since $M_0$ is irreducible, the operator $R_0$ acting on the bivectors has no nontrivial invariant
subspaces, so $\rho$ is a multiple of identity, that is, $M^n$ is Einstein, so $R(X,Y)=R_0(X,Y) + a X \wedge Y$, for some constant
$a$. Then, as $R(X,Y).R=0$ (where $R(X,Y)$ acts as a derivation), we have $R(X,Y).R_0=0$. On the other hand, $R_0(X,Y).R_0=0$, so
$a (X \wedge Y).R_0=0$. If $a \ne 0$, it would follow that the holonomy algebra of $R_0$ is the whole $\og(n)$, so $M_0$ would
have a constant curvature, which contradicts the fact that $W_0 \ne 0$. So $a=0$, hence the curvature tensors of the locally
symmetric spaces $M^n$ and $M_0$ are equal, which implies that they are isometric. 

\section{Proof of the Proposition}
\label{s:alg}

Let $\g$ be a compact simple Lie algebra with the inner product being the negative of the Killing form. Introduce the inner product
$\<A_1, A_2\>=\Tr(A_1 A_2^t)$ on $\End(\g)$.

Equation \eqref{eq:biad} is easily seen to be satisfied if $\Phi=0$ and $K_Z \in \ad(\g)$, for all $Z$. It follows
that for any solution $(K, \Phi)$ of linear system (\ref{eq:cyclePhi}, \ref{eq:biad}), $(\pi_{\ad(\g)^\perp}K, \Phi)$ is also a
solution of that system, where $\ad(\g)^\perp$ is the orthogonal complement to $\ad(\g) \subset \End(\g)$. We can therefore assume for the rest of
the proof that
\begin{equation}\label{eq:Kperpad}
K_Z \perp \ad(\g), \quad \text{for all} \quad Z \in \g,
\end{equation}
and will be proving that $K=0$.

For an orthonormal basis $\{e_i\} \subset \g$, we abbreviate $\ad_{e_i}$ to just $\ad_i$. We have
$\sum_i \ad_i^2 = -\id$ (minus the Casimir operator of $\ad$). Introduce the coboundary operator $\db:\End(\g) \to \g$ by
$\db(A) = -\frac12 \sum_i [Ae_i,e_i]$. Then 
\begin{equation}\label{eq:db}
\db(X \wedge Y) = [X, Y], \quad \db(\ad_X)= \tfrac12 X, \quad \Ker \db = \ad(\g)^\perp,
\end{equation}
for all $X, Y \in \g$, where the last two equations follow from the fact that
$\<\db(A),Y\>=\frac12\sum_i\<Ae_i,[Y,e_i]\>=-\frac12\Tr (A \ad_Y)=\frac12\<A, \ad_Y\>$, for any $A \in \End(\g)$.

We start with the following Lemma (the first two assertion of which are probably well-known).

{ 
\begin{lemma}\label{l:LPsi}
\begin{enumerate}[1.] Let $\g$ be a simple compact Lie algebra with a bi-invariant inner product.
    \item \label{it:L}
    Suppose that $L \in \End(\g)$ satisfies $\sigma_{XYZ}\<LX,[Y,Z]\>=0$, for all $X, Y, Z \in \g$. Then $L=\ad_U$,
    for some $U \in \g$.
    \item \label{it:Psi}
    Suppose that $\Xi \in \Hom(\Lambda^2\g, \g)$ satisfies $\sigma_{XYZ}[\Xi(X,Y),Z]=0$, for all $X, Y, Z \in \g$.
    Then $\Xi(X,Y)=c[X,Y]$, for some $c \in \br$.
    \item
    In the assumptions of the Proposition and condition \eqref{eq:Kperpad}, we have:
    \begin{gather}\label{eq:Phi}
    \Phi(X,Y) = \tfrac12(K_YX-K_XY), \quad \text{for all} \quad X, Y \in \g,\\
    \sum\nolimits_j \ad_j K_{e_j}=0, \label{eq:sumadK} \\
    \sigma_{XYZ}\<K_XY,Z\>=0, \quad \text{for all} \quad X, Y, Z \in \g. \label{eq:cycleK}
    \end{gather}
\end{enumerate}
\end{lemma}
\begin{proof}
1. From the assumption, it follows that for all $X, Y \in \g, \; [LX,Y]+[X,LY]=-L^t[X,Y]$, so $L$ is a
quasiderivation \cite{LL}. By \cite[Lemme~4]{Be} (or \cite[Corollary~4.14]{LL}), $L \in \ad(\g) \oplus \br \,\id$ as the
centroid of a simple Lie algebra (the centralizer of $\ad(\g)$ in $\End(\g)$) is $\br \,\id$, which then implies that $L \in \ad(\g)$.

2. From the assumption, for any $U \in \g$,
$\<[\Xi(X,Y),Z], U\>=-\<[\Xi(X,Y),U], Z\>=\<[\Xi(U,X),Y]+[\Xi(Y,U),X], Z\>
=\<[L_UX,Y]-[L_UY,X],Z\>=\<L_UX,[Y,Z]\>+\<L_UY,[Z,X]\>$, where
$L_UX=\Xi(U,X)$, so taking the inner product of the equation $\sigma_{XYZ}[\Xi(X,Y),Z]=0$ with $U$ we get 
$\sigma_{XYZ}(\<L_UX,[Y,Z]\>)=0$, for all $U \in \g$, which by assertion~\ref{it:L} implies that $L_U \in \ad(\g)$.
Then by linearity, $\Xi(U,X)=[AU,X]$, for all $X, U \in \g$, where $A \in \End(\g)$. As $\Xi$ is skew-symmetric,
we have $\ad_{AU}=\ad_UA$ (which means that $A$ is in the quasicentroid of $\g$), 
so $-A=\sum_i \ad_i^2 A=\sum_i \ad_i \ad_{Ae_i}$. Then $A^t-A=\sum_i \ad_{[e_i,Ae_i]}=0$ (as $[AX,X]=\Xi(X,X)=0$),
so $A$ is symmetric. But then from $[AU,X]+[AX,U]=0$, for any two eigenspaces $\g_\a, \g_\b$ of $A$ with the
eigenvalues $\la_\a \ne \la_\b$, we get $[\g_\a, \g_\b]=0$, so $\g_\a$ is an ideal. As $\g$ is simple, it follows that
$A=c\,\id$, for some $c \in \br$.

3. Acting by $\db$ on \eqref{eq:biad} and using \eqref{eq:Kperpad} (which implies $[\ad_{[X, Y]},K_Z] \perp \ad(\g)$)
we get from \eqref{eq:db}:
\begin{equation}\label{eq:sigmapsi}
\sigma_{XYZ}[\Xi(X,Y),Z]=0, \quad \text{where} \quad \Xi(X,Y)=K_XY-K_YX + 2 \Phi(X,Y).
\end{equation}
Then by assertion~\ref{it:Psi}, $\Xi(X,Y)=c[X,Y]$, so $\Phi(X,Y) = \frac{c}{2}[X,Y]+\frac12(K_YX-K_XY)$, which by
\eqref{eq:cyclePhi} implies $\sigma_{XYZ}(\<K_XY,Z\>)=\frac{3c}{2}\<[X,Y],Z\>$. Then for all $j=1, \dots, n$,
$K_Xe_j=K_{e_j}X-\sum_i \<K_{e_i}X,e_j\>e_i+\frac{3c}{2}[X,e_j]$. As $\db K_X=0$ by (\ref{eq:Kperpad}, \ref{eq:db}),
$0=\sum_j[K_{e_j}X-\sum_i \<K_{e_i}X,e_j\>e_i+\frac{3c}{2}[X,e_j],e_j]=2\sum_j[K_{e_j}X,e_j]-\frac{3c}{2}X$. Then
$\sum_j \ad_j K_{e_j}=-\frac{3c}{4}\,\id$. Since $\Tr (\ad_j K_{e_j})=0$ from \eqref{eq:Kperpad}, we have $c=0$, which
implies equations (\ref{eq:Phi}, \ref{eq:sumadK}) and \eqref{eq:cycleK}.
\end{proof}
} 

In the next two Lemmas we will be working with $\g^\bc$, the complexification of $\g$. We extend the maps
$K, \Phi$ and $\wedge$ by the complex linearity, with the latter one defined using the complex inner product on $\g^\bc$, the negative of the
Killing form. We will use the same notation for the complexified maps. Note that equations
(\ref{eq:cyclePhi}, \ref{eq:biad}) and \eqref{eq:Phi} still hold in $\g^\bc$.

The next Lemma proves the proposition for all the algebras of rank at least three.

{ 
\begin{lemma}\label{l:Phi0}
In the assumptions of the Proposition and condition \eqref{eq:Kperpad}, suppose that $\rk \g \ge 3$. Then
\begin{enumerate}[1.]
    \item \label{it:comm}
    $\Phi(X,Y)=0$, for all $X, Y \in \g^\bc$ such that $[X, Y]=0$.
    \item $\Phi(X,Y)=0$ and $K_XY=0$, for all $X, Y \in \g$.
\end{enumerate}
\end{lemma}
\begin{proof}
1. Let $X,Y,Z \in \g^\bc$ span a three-dimensional abelian subalgebra $\ag_3 \subset \g^\bc$.
Substituting such $X,Y,Z$ to \eqref{eq:biad} we obtain
\begin{equation}\label{eq:adsigma}
\ad_{U}=- \sigma_{XYZ}(\Phi(X,Y)\wedge Z), \quad \text{where} \quad U=\sigma_{XYZ}[K_XY-K_YX, Z].
\end{equation}
It follows that $\rk \ad_U \le 6$.

By \cite[Section~3]{Jmin}, the number $\mathbf{m}(\g^\bc)$, the minimal possible rank of $\ad_V,\; V \in \g^\bc\setminus\{0\}$ (which equals to the minimal dimension of
the nonzero adjoint orbit), is attained on the adjoint orbit of highest root vector of $\g^\bc$, if $\g^\bc \ne \mathfrak{sl}(l,\bc)$, and is attained
simultaneously on the adjoint orbit of the highest root vector and of the vectors $V=\la (-l E_{11}+ \id)$, if $\g^\bc=\mathfrak{sl}(l,\bc)$
(here $\la \in \bc, \; \la \ne 0$, and $E_{ij} \in \mathfrak{gl}(l,\bc)$ is the matrix having $1$ in the $(i,j)$-th entry and zero elsewhere).
In the both cases, $\mathbf{m}(\g^\bc)$ can be computed as in \cite[Lemma~4.3.5]{CMc}; the explicit values of $\mathbf{m}(\g^\bc)$ are given in the third row of Table~\ref{table3}
(taken from \cite[Table~1]{Jr}).

\begin{table}[h]
\renewcommand{\arraystretch}{1.2}
\begin{tabular}{|c|@{\extracolsep{.5cm}}cccccccc|}
  \hline
  $\g$ & $\mathfrak{su}(l), \; l \ge 2$ & $\mathfrak{so}(l), \; l \ge 7$ & $\mathfrak{sp}(l), \; l \ge 2$ & $\mathfrak{e}_6$ & $\mathfrak{e}_7$ & $\mathfrak{e}_8$ & $\mathfrak{f}_4$ & $\mathfrak{g}_2$ \\
  \hline
  $\rk \g$ & $l-1$ & $[l/2]$ & $l$ & $6$ & $7$ & $8$ & $4$ & $2$ \\
  $\mathbf{m}(\g^\bc)$ & $2(l-1)$ & $2(l-3)$ & $2l$ & $22$ & $34$ & $58$ & $16$ & $6$ \\
  \hline
\end{tabular}
\vspace{.2cm}
\caption{Simple compact Lie algebras, their ranks and the numbers $\mathbf{m}(\g^\bc)$.}
\label{table3} 
\end{table}


Consider two cases.

Suppose that $\rk \g \ge 4$. The inspection of Table~\ref{table3} then shows that $\mathbf{m}(\g) \ge 8$, so by
\eqref{eq:adsigma}, $\ad_U=0$. Then $\sigma_{XYZ}(\Phi(X,Y)\wedge Z)=0$, which implies $\Phi(X,Y) \in \ag_3$,
for all $X, Y, Z$ spanning a three-dimensional abelian subalgebra $\ag_3$. Taking linearly independent $X, Y$ in a
Cartan subalgebra $\h \subset \g^\bc$ and considering different $\ag_3 \subset \h$ containing $X$ and $Y$ we obtain
that $\Phi(X,Y) \in \Span(X,Y)$, for any $X,Y \in \h$. As $\Phi$ is bilinear and skew-symmetric (and the Killing form
is non-degenerate on $\h$), there exists $p \in \h$ such that $\Phi(X,Y)=\<X,p\>Y-\<Y,p\>X$. But then 
$\sigma_{XYZ}(\Phi(X,Y)\wedge Z)=0$ implies $\<X,p\>=0$, for all $X \in \h$. It follows that $\Phi(X,Y)=0$, for all $X, Y \in \h$.

Suppose that $\rk \g = 3$, so that
$\g^\bc=\mathfrak{sl}(4,\bc)(=\mathfrak{so}(6,\bc)), \; \mathfrak{so}(7,\bc)$ or $\mathfrak{sp}(3,\bc)$ (see
Table~\ref{table3}), and in all three cases, $\mathbf{m}(\g) = 6$. Acting by the both sides of \eqref{eq:adsigma} on $\g^\bc$ we get
$[U,\g^\bc] \subset \Span(\ag_3, \Phi(\ag_3,\ag_3))=\Span(X,Y,Z, \Phi(X,Y), \Phi(Y,Z), \Phi(Z,X))$. If $U \ne 0$,
the leftmost subspace has dimension at least six, so the inclusion is, in fact, an equality. In particular,
$[U,\g^\bc] \supset \ag_3$, and moreover, either $U=X_\a$, where $\a$ is the highest positive root (relative to some choice of a Cartan
subalgebra in $\g^\bc$ and positive roots) and $X_\a$ is the corresponding root vector, or $\g^\bc=\mathfrak{sl}(4,\bc)$) and $U$ lies 
in the cone of the adjoint orbit of $U_0=\diag(3,-1-1,-1)$. In the first case, the space $[U,\g^\bc]$ contains a codimension one subspace 
of nilpotent elements spanned by some positive root vectors (in fact, $[U,\g^\bc]$ is a one-dimensional solvable extension of the
five-dimensional Heisenberg algebra). It follows that if $\ag_3$ is a Cartan subalgebra, the inclusion $[U,\g^\bc] \supset \ag_3$ is
impossible. In the second case, $[U,\g^\bc]$ also contains no Cartan subalgebra $\h$ of $\g^\bc=\mathfrak{sl}(4,\bc)$ (as any $\h$ 
contains nonsingular matrices, but $[U_0,\g^\bc]$ does not).

Thus $\sigma_{XYZ}(\Phi(X,Y)\wedge Z)=0$ for all $X, Y, Z$ spanning a Cartan subalgebra $\h \subset \g^\bc$. It
follows that $\Phi(X,Y) \in \h$, for any Cartan subalgebra $\h \subset \g^\bc$ containing $X, Y$. We want to show
that $\Phi(X,Y)=0$ when $X,Y \in \h$. Consider the case $\g^\bc=\mathfrak{sl}(4,\bc)$ and take
$\h=\Span_{i=1}^4(D_i)$, where $D_i=-4E_{ii}+\id$. As every pair $D_i, D_j, \; i \ne j$, is contained in more than one Cartan 
subalgebra of $\mathfrak{sl}(4,\bc)$, we have $\Phi(D_i, D_j) \in \Span(D_i, D_j)$. As $\Phi$ is skew-symmetric, we have
$\Phi(D_i, D_j)=a_{ij}D_i-a_{ji}D_j$, for all $1 \le i \ne j \le 4$. Then $\sigma_{ijk}(\Phi(D_i,D_j)\wedge D_k)=0$
implies $a_{ki}+a_{ji}=0$, for any pairwise nonequal $i,j,k$. It follows that all the $a_{ij}$'s vanish, so
$\Phi(\h,\h)=0$. In the case $\g^\bc=\mathfrak{so}(7,\bc)$, every Cartan subalgebra $\h$ is contained in some
$\mathfrak{sl}(4,\bc)=\mathfrak{so}(6,\bc) \subset\mathfrak{so}(7,\bc)$, so from the above arguments, $\Phi(\h,\h)=0$.
In the case $\g^\bc= \mathfrak{sp}(3,\bc)=\{\bigl(\begin{smallmatrix}A&B\\C&-A^t\end{smallmatrix}\bigr) \; : \;
A,B,C \in \mathfrak{gl}(3,\bc),$ $\, B^t=B, C^t=C\}$, take $\h$ to be the diagonal subalgebra and let $D_i=E_{ii}-E_{i+3,i+3}, \; i=1,2,3$. Let
$D_4=D_1+D_2+D_3$. Then every pair $D_i, D_j, \; 1 \le i \ne j \le 4$ is contained in more than one Cartan subalgebra
of $\mathfrak{sp}(3,\bc)$, so we have $\Phi(D_i, D_j) \in \Span(D_i, D_j)$, which by the same arguments as above
implies that $\Phi(\h,\h)=0$.

So for all the algebras $\g$ of rank at least $3$, $\Phi(\h,\h)=0$, for any Cartan subalgebra $\h \subset \g^\bc$.
Then by \cite[Theorem~A]{R}, $\Phi(X,Y)=0$, for all $X, Y \in \g^\bc$ such that $[X, Y]=0$.

2. From assertion~\ref{it:comm} and by \cite[Corollary~5.2]{K}, there exists $A \in \End(\g^\bc)$ such that for all
$X, Y \in \g^\bc$, $\Phi(X,Y)=A[X,Y]$.
As $\Phi(X,Y)$ is real, when $X, Y$ are real (belong to $\g$), there exists $A \in \End(\g)$ such that for all
$X, Y \in \g, \; \Phi(X,Y)=A[X,Y]$. Then from \eqref{eq:cyclePhi} and assertion~\ref{it:L} of Lemma~\ref{l:LPsi},
$A=\ad_U$, for some $U \in \g$. Then from \eqref{eq:Phi}, $2[U,[X,e_j]] = K_{e_j}X-K_Xe_j$. Acting by $\ad_j$ and
summing up for $j=1, \dots, n$ we get $\sum_j\ad_j[U,[X,e_j]] = 0$ by \eqref{eq:sumadK} and
(\ref{eq:Kperpad}, \ref{eq:db}). Taking the inner product with an arbitrary $Y \in \g$ we obtain
$\Tr(\ad_U\ad_Y\ad_X)=0$, for all $X, Y \in \g$. It follows that $U \perp [\g, \g]=\g$, so $U=0$, hence $\Phi=0$.
Then from \eqref{eq:Phi} (and the definition of $K$) we obtain that the trilinear form $\<K_XY,Z\>$ is symmetric in
the first two variables and is skew-symmetric in the second two, so $K=0$.
\end{proof}
} 

To finish the proof of the Proposition, it remains to consider the algebras $\g$ of rank two. There are three of them
(see Table~\ref{table3}), with $\g^\bc=\mathfrak{sl}(3,\bc), \mathfrak{sp}(2,\bc)$ or $\g^\bc_2$.

{ 
\begin{lemma}\label{l:rank2}
In the assumptions of the Proposition and condition \eqref{eq:Kperpad}, suppose that $\rk \g = 2$ (so that
$\g^\bc=\mathfrak{sl}(3,\bc), \mathfrak{sp}(2,\bc)$ or $\g^\bc_2$). Then
\begin{enumerate}[1.]
    \item \label{it:KZXY}
    If $X, Y$ span a Cartan subalgebra $\h \subset \g^\bc$ and $Z \perp \h$, then $\<K_ZX,Y\>=0$.
    \item \label{it:phi}
    There exists a function $\phi: \g^\bc \to \bc$, which is homogeneous of degree one, whose restriction to
    the centralizer of any nonzero element of $\g^\bc$ is linear, and such that for all $X, Y$ with $[X,Y]=0$,
    \begin{equation}\label{eq:Phiphi}
    \Phi(X,Y)=\phi(Y)X-\phi(X)Y.
    \end{equation}
    \item \label{it:philin}
    The function $\phi$ introduced in assertion \ref{it:phi} is zero for $\g^\bc=\mathfrak{sp}(2,\bc), \g^\bc_2$
    and is linear for $\g^\bc=\mathfrak{sl}(3,\bc)$.
    \item
    $\Phi(X,Y)=0$ and $K_XY=0$, for all $X, Y \in \g$.
\end{enumerate} 
\end{lemma}
\begin{proof}
Using \eqref{eq:Phi} we can simplify \eqref{eq:biad} to the form
\begin{equation}\label{eq:biads}
\sigma_{XYZ}([\ad_{[X, Y]}+\tfrac12 X \wedge Y,K_Z]-\ad_{\db[X \wedge Y,K_Z]})=0,
\end{equation}
for all $X,Y,Z \in \g^\bc$.

1. Let $\h = \Span(U, V)\subset \g^\bc$ be a Cartan subalgebra, with $\Delta \subset \h^*$ the set of roots. For
$\a \in \Delta$, let $X_\a$ be a corresponding nonzero root vector. For $\a,\b,\gamma \in \Delta$, take
$(X,Y,Z)=(X_\a,X_\b,X_\gamma)$ in \eqref{eq:biads}.  The inner product of the resulting equation with
$U \wedge V$ gives $\sigma_{\a\b\gamma}\<K_{X_\a}[X_\b, X_\gamma], (\b+\gamma)(U)V-(\b+\gamma)(V)U\>=0$,
Choose $\a,\b,\gamma \in \Delta$ in such a way that $[X_\a,X_\b],[X_\a,X_\gamma] \in \h$ (or zero) and
$[X_\b,X_\gamma]=X_\K \ne 0$, where $\K=\b+\gamma \in \Delta$. Then we obtain
\begin{equation}\label{eq:roots}
\<K_{X_\a}X_\K, U_\K\>=0,
\end{equation}
where $U_\K$ is a nonzero vector from $\Ker \K \subset \h$.

Consider three cases separately.

For $\g^\bc=\mathfrak{sl}(3,\bc)$, for any three roots $\a_1,\a_2,\a_3$ with
$\a_1+\a_2+\a_3=0$, we can take $\a=\a_i, \; \b=-\a_i$, $\gamma=-\a_j, \; \{i, j\}=\{2,3\}$. Then
$[X_\a,X_\b]\in \h, \, [X_\a,X_\gamma]=0, \; \K=\a_k, \; \{i,j,k\}=\{1,2,3\}$, so \eqref{eq:roots} gives
$\<K_{X_\a}X_\K, U_\K\>=0$, for any $\a,\K \in \Delta$ with $\a+\K \in \Delta$. As the choice of $\h$ was arbitrary,
this equation holds under the adjoint action of $\mathrm{SL}(3,\bc)$. Fix $\K$ and act by the stabilizer of $U_\K$. 
The matrix multiplication and the linearity of $K$ show that
$\<K_XY, U_\K\>=0$, for any $Y \in \g^{U_\K}$ and for any $X \perp \g^{U_\K}$, where $\g^{U_\K}$ is the centralizer
of $U_\K$ in $\g^\bc$. Then the linearity of $K$ implies that $\<K_XV, U\>=0$, for any $U, V$ spanning a Cartan
subalgebra $\h$ and any $X \perp \h$.

For $\g^\bc=\mathfrak{sp}(2,\bc)$, the roots are $(\ve_1,\ve_2), \, (2\ve_1,0), \, (0,2\ve_2) \; i, j =\pm1$.
Choose $\K=(\ve_1,\ve_2)$. Then for any triple
$(\b,\gamma)=((2\ve_1,0),(-\ve_1,\ve_2)), \; \a=(\ve_1,-\ve_2),(-2\ve_1,0),(0,2\ve_2)$ and
$(\b,\gamma)=((0,2\ve_2),(\ve_1,-\ve_2))$, $\a=(-\ve_1,\ve_2),(2\ve_1,0),(0,-2\ve_2)$, the conditions
$[X_\a,X_\b]=[X_\a,X_\gamma] \in \h$ (or $0$) and $[X_\b,X_\gamma]=X_\K \ne 0$ are satisfied. It follows that equation
\eqref{eq:roots} is satisfied for every $\K=(\ve_1,\ve_2)$ and every $\a \in \Delta, \; \a \ne \pm\K$, so by linearity,
$\<K_XX_\K, U_\K\>=0$, for every $\K=(\ve_1,\ve_2)$ and every $X \perp \g^{U_\K}$, where $\g^{U_\K}$ is the centralizer
of $U_\K$ in $\g^\bc$. Acting by the stabilizer of $U_\K$ (and using the linearity of $K$ and the fact that
$(\g^{U_\K})^\perp$ is invariant with respect to that stabilizer) we obtain $\<K_XV, U\>=0$, for any $U, V$ spanning a
Cartan subalgebra $\h$ and any $X \in \oplus_{\a \ne \pm \K} \g_\a$, where $\g_\a$ are the root spaces and
$\K=(\ve_1,\ve_2)$ is arbitrary. By linearity it then follows that $\<K_XV, U\>=0$, for any $U, V$ spanning $\h$ and
any $X \perp \h$.

For $\g^\bc=\g^\bc_2$, the roots are $\pm \omega_i, \, \omega_i-\omega_j, \; i,j=1,2,3, \, i \ne j$, where
$\sum_i\omega_i=0$. Choose $\K=\omega_1$. Then for any triple $\a, \b,\gamma$ such that
$\b=\omega_1-\omega_i,\, \gamma=\omega_i,\; \a=-\omega_i,\omega_i-\omega_1,\omega_1-\omega_j,\; i,j=1,2,3,\, i \ne j$,
the conditions $[X_\a,X_\b]=[X_\a,X_\gamma] \in \h$ (or $0$) and $[X_\b,X_\gamma]=X_\K \ne 0$ are satisfied, so by
\eqref{eq:roots}, $\<K_XX_{\omega_1},U_{\omega_1}\>=0$, for all
$X \in S= \oplus\{\g_a, \a=-\omega_i,\pm (\omega_1-\omega_i),\; i \ne 1\}$.
The adjoint action of the subgroup $G_1=\exp(tX_{\omega_1})$ keeps $X_{\omega_1}$ and $U_{\omega_1}$ fixed, so the
equation $\<K_XX_{\omega_1},U_{\omega_1}\>=0$ holds for all
$X \in S'=\oplus\{\g_a, \a=\pm\omega_i,\pm (\omega_1-\omega_i),\; i \ne 1\}$, the smallest $G_1$-submodule of $\g^\bc$
containing $S$. As $S'$ is also an $\mathrm{SL} (2, \bc)$-submodule, where $\mathrm{SL} (2, \bc) \subset \g^\bc$ is
the subgroup tangent to $\mathfrak{sl}(2,\bc)=\Span(X_{\omega_1},X_{-\omega_1},[X_{\omega_1},X_{-\omega_1}])$, and the
adjoint action of this $\mathrm{SL} (2, \bc)$ keeps $U_{\omega_1}$ fixed, we obtain that $\<K_XY,U_{\omega_1}\>=0$,
for all $X \in S', \; Y \in \mathfrak{sl}(2,\bc)$, so, in particular, $\<K_XV, U\>=0$, for any $U, V$ spanning a
Cartan subalgebra $\h$ and any $X \in S'=\oplus\{\g_a, \a=\pm\omega_i,\pm (\omega_1-\omega_i),\; i \ne 1\}$. Repeating
the arguments with $\omega_1$ replaced by $\omega_2$ we obtain that $\<K_XV, U\>=0$, for any $U, V$ spanning $\h$ and
any $X \perp \h$.

2. From assertion \ref{it:KZXY} and by (\ref{eq:Phi}, \ref{eq:cycleK}), $\Phi(X,Y) \in \h$, if $\h=\Span(X, Y)$ is a
Cartan subalgebra. Then
$\Phi(X,Y) \wedge X \wedge Y=0$, for all $(X,Y) \in \h \times \h \subset \g^\bc \times \g^\bc$, so by \cite[Theorem~A]{R},
$\Phi(X,Y)\in \Span(X, Y)$, for any commuting $X, Y \in \g^\bc$. Denote
$\cp= \{(X,Y) \in \g^\bc \times \g^\bc \, : \, [X,Y]=0\}$ the commuting variety of $\g^\bc$ and $\cp^0$ its subset
consisting of the linearly independent pairs $(X,Y)$. As $\rk \g =2$, the subset $\cp^0$ is nonempty and is dense in $\cp$.
Since $\Phi(X,Y)\in \Span(X, Y)$, for $(X,Y) \in \cp$, there exist (uniquely defined) functions $f,g:\cp^0 \to \bc$
such that $\Phi(X,Y)=f(X,Y)X+g(X,Y)Y$, for all $(X,Y) \in \cp^0$. As $\Phi$ is skew-symmetric, $g(X,Y)=-f(Y,X)$, so
$\Phi(X,Y)=f(X,Y)X-f(Y,X)Y$. Denote $\g^X$ the centralizer of $X \ne 0$ in $\g^\bc$. If
$\dim \g^X=2, \; \g^X = \Span(X,Y)$, then from $\Phi(X,Y+tX)=\Phi(X,Y)$ it follows that $f(Y,X)=f(Y+tX,X)$, for all
$t \in \bc$. If $\dim \g^X >2$ and $Y_1, Y_2 \in \g^X$ are such that $\rk(X,Y_1,Y_2)=3$, then the linearity of $\Phi$
by $Y$ implies that $f(t_1Y_1+t_2Y_2,X)=f(Y_1,X)$. In the both cases, it follows that $f(Y,X)$ does not depend on $Y$
for $(X,Y) \in \cp^0$, so there exists a function $\phi:\g^\bc\setminus\{0\}$ such that $\Phi(X,Y)=\phi(Y)X-\phi(X)Y$,
for all $(X,Y) \in \cp^0$. The function $\phi$ is homogeneous of degree $1$ and, putting $\phi(0)=0$ we obtain that
equation \eqref{eq:Phiphi} is satisfied, for all $X, Y$ with $[X,Y]=0$. From \eqref{eq:Phiphi} it easily follows that
the restriction of the function $\phi$ to the centralizer $\g^Z$ of any nonzero $Z \in \g^\bc$ is linear.


3. For $X \in \g^\bc$, let $X=X^S+X^N$ be the Jordan decomposition, with $X^S, X^N$ the semisimple and the nilpotent parts of
$X$, respectively. As $X^S, X^N \in \g^X, \quad \phi(X)=\phi(X^S)+\phi(X^N)$. The element $X^S$ (if it is nonzero) lies in
a Cartan subalgebra $\h$, so for some two linearly independent elements
$H_\a=[X_\a, X_{-a}], \; H_\b=[X_\b, X_{-b}]$ from $\h, \; X^S=\mu_\a H_\a+\mu_\b H_\b$ and
$\phi(X^S)=\mu_\a\phi(H_\a)+\mu_\b\phi(H_\b)$, since $\phi$ is homogeneous of degree one and $H_\a, H_\b$ lie in $\g^{X^S}$.
The subalgebra $\mathfrak{sl}(2,\bc)=\Span(H_\a,X_\a, X_{-a})$ lies in the centralizer of $\Ker \a$, so the restriction
of $\phi$ to it is linear. Hence $\phi(H_\a)=\frac12\phi(H_\a+(X_\a-X_{-a}))+$ $\frac12\phi(H_\a-(X_\a-X_{-a}))$, and both
$H_\a\pm(X_\a-X_{-a})$ are nilpotent. Using similar arguments for $H_\b$ we obtain that every $X \in \g^\bc$ is a sum
of nilpotent elements $X_i$ such that $\phi(X)=\sum_i \phi(X_i)$. Now, every nilpotent element $X_i$ lies in the
nilradical of a Borel subalgebra $B_i$. As the nilradical lies in the centralizer of the highest root vector, the restriction
of $\phi$ to it is linear, so $X_i=\sum_{\a \in \Delta^+}X_{i\a}$ and $\phi(X_i)=\sum_{\a \in \Delta^+}\phi(X_{i\a})$.

To prove that $\phi=0$ for $\g^\bc=\mathfrak{sp}(2,\bc)$ and $\g^\bc=\g^\bc_2$ it is therefore sufficient to prove that $\phi(X_\a)=0$,
where $X_\a$ is a root vector for some Cartan subalgebra of $g^\bc$. By \cite{M}, every $X_\a$ can be included in a
three-dimensional abelian subalgebra $\ag \subset g^\bc$. Indeed, for $\mathfrak{sp}(2,\bc)$, with the roots
$\pm \omega_1\pm\omega_2, \pm 2 \omega_1, \pm \omega_2$, 
we can take $\ag=\Span(X_{2\ve_1\omega_1},X_{2\ve_2\omega_2},X_{\ve_1\omega_1+\ve_2\omega_2})$, where $\ve_1, \ve_2 = \pm 1$;
and for $\g^\bc_2$, with the roots $\pm \omega_i, \, \omega_i-\omega_j, \; i,j=1,2,3, \, i \ne j, \; \sum_i\omega_i=0$, we can
take $\ag=\Span(X_{\ve\omega_i},X_{\ve(\omega_i-\omega_j)},X_{\ve(\omega_i-\omega_k)})$, where $\ve = \pm 1$ and
$\{i,j,k\}=\{1,2,3\}$. Now suppose that $X,Y,Z$ span a three-dimensional abelian subalgebra of $\g^\bc$.
Substituting such $X, Y, Z$ to \eqref{eq:biads} and using the fact that
$\sigma_{XYZ}[X \wedge Y,K_Z]=\sigma_{XYZ}(-K_ZX \wedge Y+$ $K_ZY\wedge X)=2\sigma_{XYZ}(\Phi(X,Y)\wedge Z)
=4\sigma_{XYZ}(\phi(Y)X\wedge Z)$ (by (\ref{eq:Phi},\ref{eq:Phiphi})), so $\sigma_{XYZ}\db[X \wedge Y,K_Z]=0$ by \eqref{eq:db},
we obtain $\sigma_{XYZ}(\phi(Y)X\wedge Z)=0$. As the vectors $X,Y,Z$ are linearly independent, we obtain
$\phi(X)=\phi(Y)=\phi(Z)=0$. Therefore $\phi(X_\a)=0$, for any root vector of any Cartan subalgebra of $\g^\bc$, so
$\phi=0$, when $\g^\bc=\mathfrak{sp}(2,\bc)$ or $\g^\bc=\g^\bc_2$.

Let $\g^\bc=\mathfrak{sl}(3,\bc)$. In the standard representation of $\mathfrak{sl}(3,\bc)$ on $\bc^3$, the vectors
$E_{ij}, \; 1 \le i \ne j \le 3$, and $H_{ij}=E_{ii}-E_{jj}, \, i=1, j=2,3$, is a basis for $\mathfrak{sl}(3,\bc)$. For a linear
function $\phi_0:\g^\bc \to \bc$, the function $\phi'=\phi-\phi_0$ is again homogeneous of degree one and is linear
on the centralizer of every nonzero element. Choose $\phi_0$ such that $\phi'$ vanishes on all the elements of the
basis $E_{ij}, \, i \ne j, \; H_{12}, H_{13}$ (then it also vanishes on all the $H_{ij}$, as all of them belong to the
centralizer of $H_{12})$. We want to show that $\phi'=0$.
The above discussion shows that every $X \in \mathfrak{sl}(3,\bc)$ can be represented as
$X=\sum_{i,\a}X_{i\a}$, such that $\phi(X)=\sum_{i,\a}\phi(X_{i\a})$ and every $X_{i\a}$ is a root vector for some Cartan
subalgebra, so it is sufficient to show that $\phi'(X)=0$, for every $X$ which is a root vector for a Cartan subalgebra.
Every such $X$ is a rank one matrix, so $X=a \otimes b$ for some nonzero $a,b \in \bc^3$ with $\<a,b\>=0$. 
Vector $a$ can be represented as a sum of two vectors $a_1, a_2 \perp b$ having at least one zero each. Both
vectors $X_i=a_i\otimes b, \; i = 1, 2$, lie in the centralizer of $X$, so $\phi'(X)=\phi'(X_1)+\phi'(X_2)$, hence it suffices
to show that $\phi'(X)=0$, for every $X$ of rank one having a zero column. Then similar arguments applied to the rows show that
it suffices to prove that $\phi'(X)=0$, for every $X$ of rank one having a zero column and a zero row. Any such $X$ is either
of the form $X=\mu E_{ij}+\nu E_{ik}$ (or $X=\mu E_{ji}+\nu E_{ki}$), with $\{i,j,k\}=\{1,2,3\}$, or
$X= \mu\nu H_{ij} \pm (\mu^2 E_{ij}-\nu^2 E_{ji}), \; i \ne j$. In the first case, $\phi'(X)=0$, as $E_{ij}, E_{ik} \in \g^X$
and $\phi'(E_{ij})=\phi'(E_{ik})=0$, in the second case, as $X, H_{ij}, E_{ij}, E_{ji} \in \g^{H_{ik}+H_{jk}}, \; k \ne i,j$,
and $\phi'$ vanishes on each of them.

4. For $\g^\bc=\mathfrak{sp}(2,\bc), \g^\bc_2$, from assertion \ref{it:philin} and \eqref{eq:Phiphi} we have $\Phi(X,Y)=0$, for
all $X, Y \in \g^\bc$ with $[X, Y]=0$, so by \cite[Corollary~5.2]{K}, for some $A \in \End(\g^\bc), \; \Phi(X,Y)=A[X,Y]$, for
all $X, Y \in \g^\bc$. Then from (\ref{eq:Phi}, \ref{eq:cycleK}), $\<K_ZX,Y\>=2\<A[X,Y], Z\>$, so $K_Z=2\ad_{A^tZ}$, which implies
$A=0$, by \eqref{eq:Kperpad}.

For $\g^\bc=\mathfrak{sl}(3,\bc)$, from assertion \ref{it:philin} and \eqref{eq:Phiphi}, there exists $l \in \g^\bc$ such that
$\Phi(X,Y)=(X \wedge Y) l$, for all $X, Y \in \g^\bc$ with $[X, Y]=0$, so by \cite[Corollary~5.2]{K}, for some
$A \in \End(\g^\bc)$, we have $\Phi(X,Y)=A[X,Y]+(X \wedge Y) l$, for all $X, Y \in \g^\bc$. Then by (\ref{eq:Phi}, \ref{eq:cycleK}),
$\<K_ZX,Y\>=2\<A[X,Y], Z\>+2\<(X \wedge Y) l, Z\>$, so $K_Z=2\ad_{A^tZ}+2l \wedge Z$. Acting by $\db$ and using
(\ref{eq:Kperpad}, \ref{eq:db}) we obtain $A^t=-2\ad_l$. But then
$\<\Phi(X,Y),Z\>=\<A[X,Y]+(X \wedge Y) l,Z\>=-2\<[Z,[X,Y]],l\>-\<(X \wedge Y) Z,l\>$, so from \eqref{eq:cyclePhi} we get
$\sigma_{XYZ}\<(X \wedge Y) Z,l\>=2\<\sigma_{XYZ}(\<X,Z\>Y) ,l\>=0$, for all $X,Y,Z \in \g^\bc$, which implies $l=0$, hence $A=0$.
\end{proof}
}

\section{An Example} 
\label{s:red} 

In general, the claim of the Theorem is false, if the model space $M_0$ is a reducible symmetric space. A counterexample is given by 
the Riemannian spaces constructed in \cite[Theorem~4.2]{BKV}. Let $I \subset \br$ be an interval containing zero, let 
$D_{ij}(w), \; i,j=1, \dots, n-1$, be an arbitrary skew-symmetric matrix of smooth functions on $I$, let $a(w), b(w)$ be arbitrary smooth 
functions without zeros on $I$, and let $\la > 0, \, \ve =\pm 1$. Then the Riemannian manifold $(M^n, ds^2)$ with the metric tensor
\begin{equation*}
ds^2=(f_{\ve}(x_1,w) dw)^2 + \sum\nolimits_{i=1}^{n-1} (dx_i +\sum\nolimits_{j=1}^{n-1} D_{ij}(w) x_j dw)^2,
\end{equation*}
where
$f_{+}(x,w)= a(w) e^{\la x}+b(w) e^{-\la x}, \; f_{-}(x,w)= a(w) \cos(\la x)+b(w) \sin(\la x)$ is curvature homogeneous and is
modeled on the symmetric space $M^2(\kappa) \times \br^{n-2}$, where $M^2(\kappa)$ is the two-dimensional space of constant curvature
$\kappa=-\ve \la^2$. It follows that it is also Weyl homogeneous (with the same model space). However, $(M^n, ds^2)$ is not in general
conformally equivalent to its model space (or to any symmetric space at all). Indeed, assuming the Riemannian manifold
$(M^n, e^{2\sigma} ds^2)$ to be symmetric, for some smooth function $\sigma=\sigma(w, x_1, \dots, x_{n-1}), \; w=x_0$, we must have
$\n^0_l W_{hijk}=0$, where $\n^0$ is the Levi-Civita connection of $e^{2\sigma} ds^2$. A direct computation shows that for
$h=i \ne j=k, \; h,i,j,k > 1$, this implies $\db \sigma / \db x_l=0$, for all $l=0, \dots, n-1$, so $(M^n, ds^2)$ by itself has to be
symmetric. But $\n_0 R_{010i}=\kappa f_{\ve}^2 D_{1i}(w)$, for $i > 1$, which is, in general, nonzero. It follows that although $(M^n, ds^2)$
is Weyl homogeneous, it is not conformally equivalent to any symmetric space.

\end{document}